\documentclass[11pt,a4paper]{article}
\title{{\bf Remarks on spectral gaps on\\
the Riemannian path  space
} }

\author{{\bf Shizan Fang$^a$\quad Bo Wu$^b$\footnote{wubo@fudan.edu.cn}
\vspace {3mm}
\footnote{Supported in part by  Creative Research Group Fund of
the National Natural Science Foundation of
China (No. 10121101) and RFDP(20040027009).}
 }\\
{\footnotesize $^a$I.M.B, BP 47870, Universit\'e de Bourgogne,
Dijon, France}\\
{\footnotesize$^b$ Department of Mathematics, Fudan University, Shanghai, China} }

\usepackage{amsmath, amsthm, amsfonts, amssymb}
\usepackage{mathrsfs}
\newtheorem{thm}{Theorem}[section]
\newtheorem{cor}[thm]{Corollary}

\newtheorem{proposition}[thm]{Proposition}

\theoremstyle{definition}

\newcommand{\scr}[1]{\mathscr #1}

\numberwithin{equation}{section} \theoremstyle{remark}

\def\R{\mathbb R}

\def\<{\langle} \def\>{\rangle}  
   \def\Id{\text{\rm{Id}}}

\def\E{\scr E}  
\def\beg{\begin}   \def\F{\scr F}
\def\Ric{\text{\rm{Ric}}} 
\def\e{\text{\rm{e}}}    \def\ra{\rightarrow} 
 \def\L{{\mathcal L}_T^x}
\def\dis \ric{\bf ric}
      \def\Ric{\text{\rm{Ric}}} \def\ric{\text{\rm{ric}}}
  
  \def\ee{\varepsilon}

\def\H{\mathbb{H}}

\def\Ric{\text{\rm Ric}}

\let\dis=\displaystyle
\date{}

\begin{document}

\maketitle

\begin{abstract} In this paper, we will give some remarks on links between the spectral gap of the Ornstein-Uhlenbeck operator on the Riemannian path space
with lower and upper bounds of the Ricci curvature on the base manifold;
 this work was motivated by a recent work of A. Naber on the characterization of the bound of the Ricci curvature by
analysis of path spaces.
\end{abstract}

\vskip 3mm


\vskip 8mm

\section{Introduction}\label{sect1}
Let $M$ be a complete smooth Riemannian manifold of dimension $d$, and $Z$ a $C^1$-vector field on $M$. We will be concerned with the diffusion operator
$$\dis L={1\over 2}(\Delta_M -Z),$$
where $\Delta_M$ is the Beltrami-Laplace operator on $M$. Let $\nabla$ be the Levi-Civita connection and $\Ric$ the Ricci curvature tensor on $M$. We will denote
\begin{equation*}\label{eq1.1}
\Ric_Z=\Ric + \nabla Z.
\end{equation*}
It is well-known that the lower bound $K_2$ of the symmetrized $\Ric_Z^s$, that is,
\begin{equation}\label{eq1.2}
\Ric_Z^s(x)={1\over 2} \Bigl( \Ric_Z(x)+ \Ric_Z^*(x)\Bigr)\geq K_2 \, \Id,
\end{equation}
where $\Ric_Z^*$ denotes the transposed matrix of $\Ric_Z$, gives the lower bound of constants in the logarithmic Sobolev inequality
with respect to the heat measure $\rho_t(x,dy)$, associated to $L$; more precisely,

\begin{equation}\label{eq1.3}
\int_M u^2(y)\log\Bigl({u^2(y)\over ||u||_{\rho_t}^2}\Bigr)\, \rho_t(x,dy)\leq 2\,{1-e^{-K_2t}\over K_2}\int_M |\nabla u(y)|^2\, \rho_t(x,dy), \quad t>0,
\end{equation}
where $ ||u||_{\rho_t}^2=\int_M u^2(y)\, \rho_t(x,dy)$.

\vskip 2mm
Given now  a finite number of times $0<t_1<\ldots< t_N$, consider the probability measure $\nu_{t_1, \ldots, t_N}$ on $M^N$ defined by
\begin{equation}\label{eq1.4}
\int_{M^N} f\, d\nu_{t_1, \ldots, t_N} =\int_{M^N} f(y_1, \ldots, y_N)\, p_{t_1}(x,dy_1)p_{t_2-t_1}(y_1, dy_2)\cdots p_{t_N-t_{N-1}}(y_{N-1}, dy_N)
\end{equation}
where $f$ is a bounded measurable function on $M^N$.
Then with respect to the correlated metric  $|\cdot|_C$ on $TM^N$ (see definition  \eqref{eq1.85} below),
the logarithmic Sobolev inequality still holds for $\nu_{t_1, \ldots, t_N}$, that is, there is a constant $C_N>0$ such that
\begin{equation}\label{eq1.5}
\int_{M^N} f^2\log\Bigl({f^2\over ||f||_{\nu_{t_1, \ldots, t_N}}^2}\Bigr)\, d\nu_{t_1, \ldots, t_N}\leq C_N\, \int_{M^N} |\nabla f|_C^2\, d\nu_{t_1, \ldots, t_N},\quad f\in C^1(M^N).
\end{equation}

It was proved in \cite{Hs3, CHL} that under the hypothesis

\begin{equation}\label{eq1.6}
\sup_{x\in M}|||\Ric_Z(x)|||<+\infty,
\end{equation}
where $|||\cdot|||$ denotes the norm of matrices, the constant $C_N$ in \eqref{eq1.5} can be bounded, that is
\begin{equation}\label{eq1.7}
\sup_{N\geq 1} C_N <+\infty.
\end{equation}
A natural question is whether \eqref{eq1.7} still holds only under Condition \eqref{eq1.2}? In a recent work \cite{Naber},
  A. Naber proved that if the uniform bound \eqref{eq1.7} holds, then the Ricci curvature of the base manifold has an upper bound.
  It is well-known that Inequality \eqref{eq1.3} implies the lower bound \eqref{eq1.2}, therefore Condition \eqref{eq1.7} implies \eqref{eq1.6}.
  The main purpose in \cite{Naber} is to get informations on $\Ric_Z$
from the  analysis of the Riemannian path space. Let's explain briefly the context.
\vskip 2mm
Let $O(M)$ be the bundle of orthonormal frames and $\pi: O(M)\rightarrow M$ the canonical projection.
Let $H_1, \ldots, H_d$ be the canonical horizontal vector fields on $O(M)$, consider the Stratanovich stochastic differential equation (SDE) on $O(M)$:
\begin{equation}\label{SDE}
du_t(w)=\sum_{i=1}^d H_i(u_t(w))\circ dw_t^i - {1\over 2} H_Z(u_t(w))dt,\quad
u_0(w)=u_0\in \pi^{-1}(x),
\end{equation}
where $H_Z$ denotes the horizontal lift of $Z$ to $O(M)$, that is, $\pi'(u)\cdot H_Z(u)=Z(\pi(u))$. It is well-known that under Condition \eqref{eq1.2},
the life-time $\tau_x$ of the SDE \eqref{SDE} is infinite. Let

\begin{equation}\label{path}
\gamma_t(w)=\pi(u_t(w)).
\end{equation}

Then $\{\gamma_t(w); t\geq 0\}$ is a diffusion process on $M$, having $L$ as generator. The probability measure $\nu_{t_1, \ldots, t_N}$
considered in \eqref{eq1.4} is the law
of $w\rightarrow(\gamma_{t_1}(w), \ldots, \gamma_{t_N}(w))$ on $M^N$. Now consider the following path space

$$W_x^T(M)=\bigl\{\gamma : [0,T]\ra M \ \hbox{continuous},\ \gamma(0)=x\bigr\}.$$

The law $\mu_{x,T}$ on $W_x^T(M)$ of $w\ra \gamma_\cdot(w)$ is called the Wiener measure on $W_x^T(M)$.
The integration by parts formula for $\mu_{x,T}$ was first estalished in the Seminal book \cite{Bismut}, then developed in \cite{FangMalliavin, ElworthyLi};
the Cameron-Martin type quasi-invariance
of $\mu_{x,T}$ was first proved by B. Driver \cite{Dr},  completed and simplified in \cite{Hs1, Hs2, ES}. By means of Cameron-Martin, we consider the space

\begin{equation*}
\H=\Bigl\{ h: [0,T]\ra \R^d\ \hbox{absolutely continuous};\ h(0)=0, |h|_\H^2=\int_0^T |\dot h(s)|_{\R^d}^2\,ds<+\infty\Bigr\}
\end{equation*}
where the dot denotes the derivative with respect to the time $t$.
Let $F:  W_x^T(M) \ra \R$ be a cylindrical function in the form: $F(\gamma)=f(\gamma(t_1), \cdots, \gamma(t_N))$ for some $N\geq1, 0\leq t_1<t_2<\cdots<t_N\leq1$, and $f\in
C_b^1(M^N)$. The usual gradient of $F$ in Malliavin calculus is defined by
\begin{equation}\label{eq1.8}
D_\tau F(\gamma(w))=\sum_{j=1}^N u_{t_j}(w)^{-1}(\partial_j f)(\gamma_{t_1}(w), \cdots, \gamma_{t_N}(w))\, {\bf 1}_{(\tau\leq t_j)},
\end{equation}
where $\partial_j$ is the gradient with respect to the $j$-th component.  The correlated norm of $\nabla f$ is

\begin{equation}\label{eq1.85}
|\nabla f|_C^2=\sum_{j,k=1}^N \langle u_{t_j}(w)^{-1}(\partial_j f), u_{t_k}(w)^{-1}(\partial_k f)\rangle\, t_j\wedge t_k,
\end{equation}

where $ t_j\wedge t_k$ denotes the minimum between $t_j$ and $t_k$. Notice that the norm $|\nabla f|_C$ is random.
The generator $\L$ associated to the Dirichlet form

$$\E(F,F)=\int_{W_x^T(M)}\Bigl(\int_0^T |D_\tau F|^2(\gamma)\, d\tau\Bigr)\, d\mu_{x,T}(\gamma)$$

 is called the Ornstein-Uhlenbeck operator.
The powerful tool of $\Gamma_2$ of Bakry and Emery \cite{BakryEmery} is not applicable to $\L$, the reason for this is
the geometry of $W_x^T(M)$ inherted from $\H$ is quite complicated, the associated ``Ricci tensor'' being a divergent object (see \cite{CM,CF, ElworthyYan}).
When the base manifold $M$ is compact, the existence of the spectral gap for $\L$ has been proved in \cite{Fang1}.
  The logarithmic Sobolev inequality for $D_\tau F$
defined in \eqref{eq1.8} has been established in \cite{AidaElworthy}, as well as in \cite{Hs3} or \cite{CHL} where the constant was estimated using the
bound of Ricci curvature tensor of the base manifold $M$. The method used in \cite{Fang1} is the martingale representation, which takes advantage the It\^o filtration;
this method has been developed in \cite{ElworthyYan} to deal with the problem of vanishing of harmonic forms on $W_x^T(M)$.
The purpose in \cite{Naber} is to proceed in the opposite direction,  to get the bound for Ricci curvature tensor of the
base manifold $M$ from the analysis of the path space $W_x^T(M)$.

\vskip 2mm
The organization of the paper is as follows. In section \ref{sect2}, we will recall briefly basic objets in Analysis of $W_x^T(M)$. On the path space $W_x^T(M)$,
there exist two type of gradients: the usual one is more related to the geometry of the base manifold, while the damped one is easy to be handled.
In section \ref{sect3}, we will make estimation of the spectral gap of $\L$ as explicitly as possible in function of
lower bound $K_2$ and upper bound $K_1$ of $\Ric$. In section \ref{sect4}, we will study the behaviour of the spectral gap $Spect(\L)$ as $T\ra 0$. Roughly speaking,
we will get the following result:
\begin{equation*}
1-{K_1T\over 2} + o(T)\leq Spect(\L)\leq 1+{K_2T\over 2} + o(T),\quad\hbox{\rm as } T\ra 0
\end{equation*}
under the following condition (\ref{eq2.0}).

\section{Framework of the Riemannian path space}\label{sect2}

We shall keep the notations of Section \ref{sect1}, and throughout this section,  $u_t(w)$ denotes always the solution of \eqref{SDE}
and $\gamma_t(w)$ the path defined in \eqref{path}.
For any $h\in \H$, we introduce first the usual gradient on the path space $W_x^T(M)$, which gives
Formula \eqref{eq1.8} when the functional $F$ is a cylindrical function. To this end, let

\begin{equation}\label{A1}
q(t,h)=\int_0^t\Omega_{u_s(w)}\Bigl(h(s), \circ dw(s)-{1\over 2} u_s(w)^{-1}Z_{\gamma_s(w)}ds\Bigr)
\end{equation}
where $\Omega_u$ is the equivariant representation of the curvature tensor on $M$. Let $\ric_Z$ be the equivariant representation of $\Ric_Z$, that is,
$$\dis \ric_Z(u)=u^{-1}\circ\Ric_Z(\pi(u))\circ u,\quad u\in O(M).$$
Consider $\hat h(w)\in\H$ defined by
\begin{equation}\label{A2}
\dot{\hat h}_t(w)=\dot h(t) +{1\over 2} \ric_Z(u_t(w))\, h(t).
\end{equation}

Let $F: W_x^T(M)\ra \R$ be a functional, we denote $\tilde F(w)=F(\gamma_\cdot(w))$. Then according to \cite{FangMalliavin}, we define

\begin{equation}\label{A3}
(D_hF)(\gamma_\cdot(w))=\Bigl\{{d\over d\ee}\tilde F\Bigl(\int_0^\cdot e^{\ee q(s,h)}\, dw(s)+\ee\, \hat h\Bigr)\Bigr\}_{\ee=0}.
\end{equation}

By \cite{Bismut, FangMalliavin}, if $F$ is a cylindrical function on $W_x^T(M)$, then
$$\dis (D_hF)(\gamma_\cdot(w))=\int_0^T \langle D_\tau F(\gamma_\cdot(w)),\dot h(\tau)\rangle\, d\tau$$
where $D_\tau F$ was given in \eqref{eq1.8}.  Consider the following resolvent equation
\begin{equation}\label{eq2.1}
{dQ_{t,s}\over dt}=-{1\over 2}\ric_Z(u_t(w))\,Q_{t,s},\quad t\geq s,\ Q_{s,s}=\Id.
\end{equation}
For a cylindrical function $F$ on $W_x^T(M)$ given by $F(\gamma)=f(\gamma(t_1),\cdots, \gamma(t_N))$ with $f\in C_b^1(M^N)$,
following \cite{FangMalliavin}, we define the damped gradient $\tilde D_\tau F$ of $F$  by
\begin{equation}\label{eq2.2}
\tilde D_\tau F(\gamma_\cdot(w))=\sum_{j=1}^N Q_{t_j,\tau}^*\bigl(u_{t_j}(w)^{-1}\partial_j f\bigr)\, {\bf 1}_{(\tau\leq t_j)},
\end{equation}
where $Q_{\tau,s}^*$ is the transpose matrix of $Q_{\tau,s}$. The damped gradient $\tilde D_\tau F$ on the path space $W_x^T(M)$ plays a basic role in Analysis of $W_x^T(M)$.
 Let $(v_t)_{t\geq 0}$ be a $\R^d$-valued process, adapted to the It\^o filtration $\F_t$
generated by $\{w(s);\ s\leq t\}$ such that ${\mathbb E}(\int_0^T |v_t|^2\, dt)<+\infty$. Consider two maps $v\ra \tilde v$ and $v\ra \hat v$ defined respectively
by
\begin{equation}\label{eq2.3}
\tilde v_t=v_t-{1\over 2}\ric_{u_t(w)}\,\int_0^t Q_{t,s}v_s\, ds,
\end{equation}
and
\begin{equation}\label{eq2.4}
\hat v_t = v_t +{1\over 2}\ric_{u_t(w)}\, \int_0^t v_s\,ds.
\end{equation}
Then $\hat{\tilde v}=\tilde{\hat v}=v$. The two gradients $D_tF$ and $\tilde D_tF$ are linked by the following formula
\begin{equation}\label{eq2.5}
\int_0^T \langle \tilde D_tF, v_t\rangle \, dt=\int_0^T\langle D_tF, \tilde v_t\rangle\ dt.
\end{equation}
The good feather of the damped gradient is that it admits a nice martingale representation
$$F=\mathbb{E}(F)+\int^T_0\<\mathbb{E}^{\F_t}(\tilde{D}_tF), d w_t\>$$
where $\mathbb{E}^{\F_t}$ denotes the conditional expectation with respect to $\F_t$. The following logarithmic Sobolev inequality holds
(\cite{ELJL, FWW}):

\begin{equation}\label{eq2.7}
\mathbb{E}\bigg(F^2\log \frac{F^2}{\|F\|^2_{L^2}}\bigg)\leq 2\mathbb{E}\bigg(\int^T_0|\tilde{D}_tF|^2dt\bigg).
\end{equation}

\section{Precise lower bound on the spectral gap}\label{sect3}

The inconvenient of Inequality \eqref{eq2.7} is that the geometric information of the base manifold $M$ is completely hidden. Now we use the usual gradient $D_tF$ to
make involving the geometry of $M$.
By \eqref{eq2.7}, the matter is now to estimate $\int^T_0|\tilde{D}_tF|^2dt$ by  $|D_tF|$.  We assume that
\begin{equation}\label{eq2.0}
K_2\,\Id\leq \dis \ric_Z^s,\quad |||\ric_Z|||\leq K_1
\end{equation}
for two constants $K_1,K_2$ with $K_1\geq0$ and $K_1+K_2\geq0$.

\begin{thm}\label{th1} Let $0<t \leq T$. Set
\begin{equation}\label{th1.1}
\begin{split}
 \Lambda(t,T)
&=1+\frac{K_1}{K_2}\Big(1-\e^{-\frac{K_2(T-t)}{2}}\Big)+\frac{K_1}{K_2}\Big(1-\e^{-\frac{K_2t}{2}}\Big)\\
&+\Big(\frac{K_1}{K_2}\Big)^2\bigg[\Big(1-\e^{-\frac{K_2t}{2}}\Big)+\frac{1}{2}\Big(\e^{-\frac{K_2(T+t)}{2}}
-\e^{-\frac{K_2(T-t)}{2}}\Big)\bigg].
\end{split}
\end{equation}
Then we have the relation:
\begin{equation}\label{th1.2}
\int^T_0|\tilde{D}_tF|^2 dt\leq \int^T_0\Lambda(t,T)|D_tF|^2dt.
\end{equation}
\end{thm}

\vskip 2mm
{\bf Proof.}
From (\ref{eq2.2}) and  (\ref{eq2.5}), we have
\begin{equation}\label{eq2.8}\tilde{D}_tF=D_tF-\frac{1}{2}\int_t^TQ_{s,t}^*\,ric^*_{u_s}D_sF ds.\end{equation}
Thus,
\begin{equation*}\label{c2.7}\aligned
|\tilde{D}_tF|^2&=|D_tF|^2-\bigg\<D_tF,\int_t^TQ_{s,t}^*ric^*_{u_s}D_sF ds\bigg\>+\frac{1}{4}\bigg|\int_t^TQ_{s,t}^*ric^*_{u_s}D_sF ds\bigg|^2\\
&:=I_1+I_2+I_3\ \hbox{\rm respectively}.\endaligned
\end{equation*}
In the following we will estimate the term of $I_2$ and $I_3$. Under the lower bound in \eqref{eq2.0},
$$\dis |||Q_{s,t}^*|||\leq e^{-{K_2(s-t)\over 2}},\quad s\geq t.$$

Let $$\Lambda_1(t,T):=\int^T_t\e^{-\frac{K_2(s-t)}{2}}ds.$$ Then

\begin{equation*}\aligned
|I_2|&\leq |D_tF|\int_t^T\e^{-\frac{K_2(s-t)}{2}}K_1|D_sF|ds\\
&\leq|D_tF|\sqrt{K_1\int^T_t\bigg(\e^{-\frac{K_2(s-t)}{4}}\bigg)^2ds}\,\sqrt{K_1\int_t^T\e^{-\frac{K_2(s-t)}{2}}|D_sF|^2ds}\\
&=|D_tF|\sqrt{K_1\Lambda_1(t,T)}\,\sqrt{K_1\int_t^T\e^{-\frac{K_2(s-t)}{2}}|D_sF|^2ds}
\\&\leq\frac{1}{2}\bigg\{|D_tF|^2K_1\Lambda_1(t,T)+K_1\int_t^T\e^{-\frac{K_2(s-t)}{2}}|D_sF|^2ds\bigg\}.
\endaligned\end{equation*}

and
\begin{equation*}\aligned
|I_3|&\leq \frac{1}{4}\bigg|\int_t^T\e^{-\frac{K_2(s-t)}{2}}K_1|D_sF|ds\bigg|^2
\\&\leq\frac{1}{4}K_1^2\Lambda_1(t,T)\int_t^T\e^{-\frac{K_2(s-t)}{2}}|D_sF|^2ds.
\endaligned\end{equation*}

Combining all the above inequalities, we get
\begin{equation*}\aligned|\tilde{D}_tF|^2&\leq \Big(1+\frac{K_1}{2}\Lambda_1(t,T)\Big)|D_tF|^2
+\Big(1+\frac{K_1}{2}\Lambda_1(t,T)\Big)\frac{K_1}{2}\int_t^T\e^{-\frac{K_2(s-t)}{2}}|D_sF|^2ds\\
&=\Big(1+\frac{K_1}{2}\Lambda_1(t,T)\Big)\bigg(|D_tF|^2+\frac{K_1}{2}\int_t^T\e^{-\frac{K_2(s-t)}{2}}|D_sF|^2ds\bigg).\endaligned\end{equation*}

Therefore, we obtain

\begin{equation*}\aligned\int^T_0|\tilde{D}_tF|^2dt&\leq\int^T_0\Big(1+\frac{K_1}{2}\Lambda_1(t,T)\Big)|D_tF|^2dt\\
&+\int^T_0\Big(1+\frac{K_1}{2}\Lambda_1(t,T)\Big)\frac{K_1}{2}\int_t^T\e^{-\frac{K_2(s-t)}{2}}|D_sF|^2dsdt\\
&=\int^T_0\Big(1+\frac{K_1}{2}\Lambda_1(s,T)\Big)|D_sF|^2ds\\
&+\int^T_0|D_sF|^2ds\int_0^s\frac{K_1}{2}\Big(1+\frac{K_1}{2}\Lambda_1(t,T)\Big)\e^{-\frac{K_2(s-t)}{2}}dt\\
&:=\int^T_0\Big(1+\frac{K_1}{2}\Lambda_1(s,T)\Big)|D_sF|^2ds+\int^T_0(J_1(s)+J_2(s))|D_sF|^2ds,\endaligned
\end{equation*}
where
$$J_1(s):=\int_0^s\frac{K_1}{2}\e^{-\frac{K_2(s-t)}{2}}dt,~J_2(s):=\int_0^s\Big(\frac{K_1}{2}\Big)^2\Lambda_1(t,T)\e^{-\frac{K_2(s-t)}{2}}dt.$$

Next, then we compute the term $J_1(s)$ and $J_2(s)$. By direct computation, we have
\begin{equation*}J_1(s)=\frac{K_1}{K_2}\Big(1-\e^{-\frac{K_2s}{2}}\Big)
\end{equation*}
and
\begin{equation*}\aligned J_2(s)&=\Big(\frac{K_1}{2}\Big)^2\int_0^s\frac{2}{K_2}
\Big(1-\e^{-\frac{K_2(T-t)}{2}}\Big)\e^{-\frac{K_2(s-t)}{2}}dt\\
&=\Big(\frac{K_1}{2}\Big)^2\frac{2}{K_2}\bigg[\int_0^s\e^{-\frac{K_2(s-t)}{2}} dt-\e^{-\frac{K_2(T+s)}{2}}\int^s_0
\e^{K_2t}dt\bigg]\\
&=\Big(\frac{K_1}{2}\Big)^2\frac{2}{K_2}\bigg[\frac{2}{K_2}\Big(1-\e^{-\frac{K_2s}{2}}\Big)-\frac{1}{K_2}\e^{-\frac{K_2(T+s)}{2}}(\e^{K_2s}-1)\bigg]\\
&=\Big(\frac{K_1}{2}\Big)^2\frac{2}{K_2}\bigg[\frac{2}{K_2}\Big(1-\e^{-\frac{K_2s}{2}}\Big)+\frac{1}{K_2}\e^{-\frac{K_2(T+s)}{2}}
-\frac{1}{K_2}\e^{-\frac{K_2(T-s)}{2}}\bigg].\endaligned
\end{equation*}
Adding $J_1(s)$ to $J_1(s)$ implying that
$$\aligned &J_1(s)+J_2(s)\\
&=\frac{K_1}{K_2}\Big(1-\e^{-\frac{K_2s}{2}}\Big)+
\Big(\frac{K_1}{K_2}\Big)^2\bigg[\Big(1-\e^{-\frac{K_2s}{2}}\Big)+\frac{1}{2}\Big(\e^{-\frac{K_2(T+s)}{2}}
-\e^{-\frac{K_2(T-s)}{2}}\Big)\bigg]:=\Lambda_2(s,T)\endaligned$$
Thus,
\begin{equation*}\int^T_0|\tilde{D}_tF|^2 dt\leq \int^T_0\Lambda(t,T)|D_tF|^2dt,\end{equation*}
with
$$\aligned \Lambda(t,T)&=1+\frac{K_1}{2}\Lambda_1(t,T)+\Lambda_2(t,T)\\
&=1+\frac{K_1}{K_2}\Big(1-\e^{-\frac{K_2(T-t)}{2}}\Big)+\frac{K_1}{K_2}\Big(1-\e^{-\frac{K_2t}{2}}\Big)\\
&+\Big(\frac{K_1}{K_2}\Big)^2\bigg[\Big(1-\e^{-\frac{K_2t}{2}}\Big)+\frac{1}{2}\Big(\e^{-\frac{K_2(T+t)}{2}}
-\e^{-\frac{K_2(T-t)}{2}}\Big)\bigg].\endaligned$$
The proof is completed. $\square$

\vskip 2mm
Notice that as $K_2\ra 0$, by expression \eqref{th1.1},
$$\Lambda(t,T)\ra 1+ {K_1T\over 2} + K_1^2 \Bigl({Tt \over 4}-{t^2\over 8}\Bigr).$$

Now we study the variation of the function $t\ra \Lambda(t,T)$. It is quite interesting to remark that its monotonicity is dependent of the sign of $K_2$.

\begin{proposition}\label{prop1}  (i) If $K_2<0$, then $t\ra \Lambda(t,T)$ is strictly increasing over $[0,T]$.
(ii) If $K_2>0$, then the maximum is attained at a point $t_0$ in $(0, T)$.
\end{proposition}

\vskip 2mm
{\bf Proof.}
Taking the derivative of $t\ra \Lambda(t,T)$ gives
$$\aligned \Lambda'(t,T)
&=-\frac{K_1}{2}\e^{-\frac{K_2(T-t)}{2}}+\frac{K_1}{2}\e^{-\frac{K_2t}{2}}\\&
+\frac{K_1^2}{2K_2}\e^{-\frac{K_2t}{2}}-\frac{K_1^2}{4K_2}\e^{-\frac{K_2(T+t)}{2}}-\frac{K_1^2}{4K_2}\e^{-\frac{K_2(T-t)}{2}}.\endaligned$$

In addition, we have
$$ \Lambda(0,T)=1+\frac{K_1}{K_2}\Big(1-\e^{-\frac{K_2T}{2}}\Big)$$
and
$$\aligned\Lambda(T,T)&=1+\frac{K_1}{K_2}\Big(1-\e^{-\frac{K_2T}{2}}\Big)+\Big(\frac{K_1}{K_2}\Big)^2\bigg[\Big(1-\e^{-\frac{K_2T}{2}}\Big)
+\frac{1}{2}\Big(\e^{-K_2T}-1\Big)\bigg]\\
&=1+\frac{K_1}{K_2}\Big(1-\e^{-\frac{K_2T}{2}}\Big)+\frac{1}{2}\Big(\frac{K_1}{K_2}\Big)^2\Big(1-\e^{-\frac{K_2T}{2}}\Big)^2\\
&=\frac{1}{2}+\frac{1}{2}\bigg[1+\frac{K_1}{K_2}\Big(1-\e^{-\frac{K_2T}{2}}\Big)\bigg]^2=\frac{1}{2}+\frac{1}{2}\Lambda^2(0,T).
 \endaligned$$

From the second equality in the above, we observe that $\Lambda(T,T)\geq\Lambda(0,T)$.
Moreover,
\begin{equation}\label{c2.9}\aligned \Lambda'(0,T)&=-\frac{K_1}{2}\e^{-\frac{K_2T}{2}}+\frac{K_1}{2}
+\frac{K_1^2}{2K_2}-\frac{K_1^2}{4K_2}\e^{-\frac{K_2T}{2}}-\frac{K_1^2}{4K_2}\e^{-\frac{K_2T}{2}}\\
&=\frac{K_1}{2}\Big(1-\e^{-\frac{K_2T}{2}}\Big)+\frac{K_1^2}{2K_2}\Big(1-\e^{-\frac{K_2T}{2}}\Big)\\
&=\frac{K_1}{2}(K_1+K_2)\frac{1-\e^{-\frac{K_2T}{2}}}{K_2}\geq0;
 \endaligned\end{equation}
and
\begin{equation}\label{c3.0}\aligned \Lambda'(T,T)&=-\frac{K_1}{2}+\frac{K_1}{2}\e^{-\frac{K_2T}{2}}
+\frac{K_1^2}{2K_2}\e^{-\frac{K_2T}{2}}-\frac{K_1^2}{4K_2}\e^{-K_2T}-\frac{K_1^2}{4K_2}\\
&=-\frac{K_1}{2}+\frac{K_1}{2}\e^{-\frac{K_2T}{2}}
-\frac{K_1^2}{4K_2}\Big(1-2\e^{-\frac{K_2T}{2}}+\e^{-K_2T}\Big)\\
&=-\frac{K_1}{2}\Big(1-\e^{-\frac{K_2T}{2}}\Big)
-\frac{K_1^2}{4K_2}\Big(1-\e^{-\frac{K_2T}{2}}\Big)^2.
 \endaligned\end{equation}
We see that

\begin{equation}\label{sign}
 \left\{ \begin{array}{ll}
\Lambda'(T,T)>0  \qquad   &if ~K_2<0,\\
\Lambda'(T,T)<0 & if ~K_2>0.
\end{array}\right.
\end{equation}

Now we look for $t\in [0,T]$ such that $\Lambda'(t,T)=0$. We have

\begin{equation}\label{c3.1}\aligned &\Lambda'(t,T)=0 \\
&\Leftrightarrow~~\Big(-\frac{K_1}{2}\e^{-\frac{K_2T}{2}}-\frac{K_1^2}{4K_2}\e^{-\frac{K_2T}{2}}\Big)\e^{K_2t}
+\Big(\frac{K_1}{2}+\frac{K_1^2}{2K_2}-\frac{K_1^2}{4K_2}\e^{-\frac{K_2T}{2}}\Big)=0\\
&\Leftrightarrow~~-\frac{K_1}{4}\e^{-\frac{K_2T}{2}}\Big(2+\frac{K_1}{K_2}\Big)\e^{K_2t}
+\frac{K_1}{4}\Big(2+\frac{2K_1}{K_2}-\frac{K_1}{K_2}\e^{-\frac{K_2T}{2}}\Big)=0
\\
&\Leftrightarrow~~\e^{-\frac{K_2T}{2}}\Big(2+\frac{K_1}{K_2}\Big)\e^{K_2t}
=\Big(2+\frac{2K_1}{K_2}-\frac{K_1}{K_2}\e^{-\frac{K_2T}{2}}\Big).\endaligned\end{equation}
Therefore there exists at most one $t$ such that  $\Lambda'(t,T)=0$.
For the case where $K_2<0$, if there exists $t_0\in(0,T)$ such that $\Lambda(t_0,T)<0$.
Then by (\ref{c2.9}) and \eqref{sign}, the equation $\Lambda'(t,T)=0$ has at least two solutions, it is impossible.
Therefore for $K_2<0$, $\Lambda'(t,T)\geq0$. For $K_2>0$, we suppose $t_0$ such that $\Lambda'(t_0,T)=0$.
Let $\beta=\frac{K_1}{K_2}$, then by (\ref{c3.1})
$$\e^{K_2t_0}=\Big(1+\frac{\beta}{2+\beta}\Big(1-\e^{-\frac{K_2T}{2}}\Big)\Big)\e^{\frac{K_2T}{2}},$$
or $t_0\in (0,T)$ is such that

\begin{equation}\label{max}
\e^{\frac{K_2t_0}{2}}=\sqrt{1+\frac{\beta}{2+\beta}\Big(1-\e^{-\frac{K_2T}{2}}\Big)}\e^{\frac{K_2T}{4}}.
\end{equation}

The proof is completed. $\square$

\begin{proposition}\label{prop2} Let $\dis \beta={K_1\over K_2}$, then (i) if $K_2>0$,
\begin{equation}\label{prop2.1}
\begin{split}
 \sup_{t\in [0,T]}\Lambda(t,T)
&=(1+\beta)^2-\Big(\beta+\frac{\beta^2}{2}\Big)\sqrt{1+\frac{\beta}{2+\beta}\Big(1-\e^{-\frac{K_2T}{2}}\Big)}\ \e^{-\frac{K_2T}{4}}\\
&-\frac{\Big(\beta+\beta^2-\frac{\beta^2}{2}\e^{-\frac{K_2T}{2}}\Big)}{\sqrt{1+\frac{\beta}{2+\beta}\Big(1-\e^{-\frac{K_2T}{2}}\Big)}}
\e^{-\frac{K_2T}{4}}.
\end{split}
\end{equation}
(ii) if $K_2<0$,
\begin{equation}\label{prop2.2}
\sup_{t\in[0,T]}\Lambda(t,T)=\frac{1}{2}+\frac{1}{2}\Big(1+\frac{K_1}{K_2}\Big[1-\e^{-\frac{K_2T}{2}}\Big]\Big)^2.
\end{equation}

\end{proposition}

\vskip 2mm
{\bf Proof.}  For $K_2>0$,
 we have
$$\aligned\Lambda(t_0,T)&=1+\beta\Big(1-\e^{-\frac{K_2T}{2}}\cdot\e^{\frac{K_2t_0}{2}}\Big)+\beta\Big(1-\e^{\frac{K_2t_0}{2}}\Big)\\
&+\beta^2\bigg[\Big(1-\e^{-\frac{K_2t_0}{2}}\Big)+\frac{1}{2}\Big(\e^{-\frac{K_2T}{2}}\cdot
\e^{\frac{-K_2t_0}{2}}-\e^{-\frac{K_2T}{2}}\cdot\e^{\frac{K_2t_0}{2}}\Big)\bigg]\\
&=1+2\beta+\beta^2-\Big(\beta+\frac{\beta^2}{2}\Big)\e^{-\frac{K_2T}{2}}\cdot\e^{\frac{K_2t_0}{2}}
-\Big(\beta+\beta^2-\frac{\beta^2}{2}\e^{-\frac{K_2T}{2}}\Big)\e^{-\frac{K_2t_0}{2}}.\endaligned$$
Using \eqref{max} yields \eqref{prop2.1}.
For $K_2<0$, $\sup_{t\in[0,T]}\Lambda(t,T)=\Lambda(T,T)$, which gives \eqref{prop2.2}. $\square$

\vskip 2mm
Combining \eqref{eq2.7} and \eqref{th1.2}, we get
\begin{thm}\label{th2} Let $\dis C(T,K_1,K_2)=\sup_{t\in [0,T]}\Lambda(t,T)$; then it holds

\begin{equation}\label{th2.1}
\mathbb{E}\bigg(F^2\log \frac{F^2}{\|F\|^2_{L^2}}\bigg)\leq 2 C(T,K_1,K_2)\mathbb{E}\bigg(\int^T_0|D_tF|^2dt\bigg)
\end{equation}
for any cylindrical function $F$ on $W_x^T(M)$.
\end{thm}

It is well-konwn that the above logarithmic Sobolev inequality implies that the spectral gap of $\L$, denoted by $Spect(\L)$, has the following lower bound
$$\dis Spect(\L)\geq {1\over C(T,K_1, K_2)}.$$

\begin{thm} Assume (\ref{eq2.0}) holds, then (i)  if $K_2>0$, we have
\begin{equation}\label{c3.2}
Spect(\L)^{-1}\leq \Bigl(1+{K_1\over K_2}\Bigr)^2-{K_1\over K_2}\sqrt{\Bigl(2+{K_1\over K_2}\Bigr)\Bigl(2+2{K_1\over K_2}-{K_1\over K_2}\e^{-\frac{K_2T}{2}}\Bigr)}\ \e^{-\frac{K_2T}{4}};
\end{equation}
(ii) if $K_2<0$, we have
\begin{equation}\label{c3.22}
Spect(\L)^{-1} \leq \frac{1}{2}+\frac{1}{2}\Big(1+\frac{K_1}{K_2}\Big[1-\e^{-\frac{K_2T}{2}}\Big]\Big)^2.
\end{equation}
\end{thm}

\vskip 2mm
{\bf Proof.}  Using the elementary inequality: $A+B\geq 2\sqrt{AB}$  to the last two terms in \eqref{prop2.1} yields \eqref{c3.2}.
Inequality \eqref{c3.22} is obvious. $\square$

\vskip 2mm
It is quite interesting to remark that

\begin{proposition}\label{prop3}
Let $\psi(T, K_1, K_2)$ be the right hand side of \eqref{c3.2} when $K_2>0$ and the right hand side of \eqref{c3.22} for $K_2<0$, then
\begin{equation}\label{prop3.1}
\psi(T, K_1, K_2)\ra 1 + {K_1T\over 2} + {K_1^2T^2\over 8}\quad\hbox{as}\quad K_2\ra 0.
\end{equation}
\end{proposition}

\vskip 2mm
{\bf Proof.} It is easy to see that the right hand side of \eqref{c3.22} tends to $1 + {K_1T\over 2} + {K_1^2T^2\over 8}$ as $K_2\ra 0$.
For the right hand side of \eqref{c3.2}, we first remark that

$$\dis {K_1\over K_2} e^{-{K_2T\over 4}}= {K_1\over K_2} -{K_1T\over 4} + {K_1K_2T^2\over 32} + o(K_2). \leqno(a)$$

Secondly
$$\aligned
&\Bigl(2+{K_1\over K_2}\Bigr)\Bigl(2+2{K_1\over K_2}-{K_1\over K_2}\e^{-\frac{K_2T}{2}}\Bigr)\\
&= \Bigl(2+{K_1\over K_2}\Bigr)\Bigl(2+{K_1\over K_2}+{K_1T\over 2}-{K_1K_2T^2\over 8} + o(K_2)\Bigr)\\
&= \Bigl(2+{K_1\over K_2}\Bigr)^2 \biggl( 1+ {{K_1T\over 2}-{K_1K_2T^2\over 8} + o(K_2)\over 2+{K_1\over K_2}}\biggr)
\endaligned$$

Therefore
$$\aligned
&\sqrt{\Bigl(2+{K_1\over K_2}\Bigr)\Bigl(2+2{K_1\over K_2}-{K_1\over K_2}\e^{-\frac{K_2T}{2}}\Bigr)}\\
&=\Bigl(2+{K_1\over K_2}\Bigr)\biggl( 1+ {1\over 2} {{K_1T\over 2}-{K_1K_2T^2\over 8} + o(K_2)\over 2+{K_1\over K_2}}-{K_1K_2T^2\over 32} + o(K_2^2)\biggr)\\
&=\Bigl(2+{K_1\over K_2}\Bigr)+ {K_1T\over 4}-{3K_1K_2T^2\over 32} + o(K_2).
\endaligned$$

Combining this with $(a)$, we get
$$\aligned
&{K_1\over K_2} e^{-{K_2T\over 4}}\sqrt{\Bigl(2+{K_1\over K_2}\Bigr)\Bigl(2+2{K_1\over K_2}-{K_1\over K_2}\e^{-\frac{K_2T}{2}}\Bigr)}\\
&=\Bigl(2+{K_1\over K_2}\Bigr){K_1\over K_2} - {K_1T\over 2} -{K_1^2T^2\over 8} + o(K_2).
\endaligned$$
Then \eqref{prop3.1} follows from the right hand side of \eqref{c3.2}. $\square$

\begin{cor}\label{cor2.1}  Assume (\ref{eq2.0}) holds.

(1) If $K_1=K_2=K>0$, then
$$\psi(T, K, K)= 4-\sqrt{3\Big(4-\e^{-\frac{KT}{2}}\Big)}\e^{-\frac{KT}{4}}\rightarrow 1~~as~~K\rightarrow 0.$$

(2) If $K_2=-K_1=-K$, then
$$\psi(T, K, -K)= {1\over 2}(1+e^{KT}).$$
\end{cor}

\vskip 2mm
{\bf Remark.} Our results improve estimates obtained in \cite{Aida}.

\section{Behaviour of $Spect(\L)$ as $T\ra 0$}\label{sect4}

In this section, we consider the case where $Z=0$. Then Condition \eqref{eq2.0} can be readed as

\begin{equation}\label{eq2.02}
K_2\,\Id\leq \ric\leq K_1\, \Id,\quad\hbox{\rm with } K_1+K_2\geq 0
\end{equation}
and SDE \eqref{SDE} is reduced to

\begin{equation}\label{SDE2}
du_t(w)=\sum_{i=1}^d H_i(u_u(w))\circ dw_t^i ,\quad
u_0(w)=u_0\in \pi^{-1}(x).
\end{equation}
The path $\gamma_t(w)=\pi(u_t(w))$ is called Brownian motion path on $M$.
Let $\rho(x,y)$ be the Riemannian distance. By \cite[p. 199]{Stroock}, there is $\ee>0$ such that
\begin{equation}\label{eq4.0}
\sup_{t\in [0,T]}{\mathbb E}\biggl(\exp\Big(\ee {\rho(x, \gamma_t)^2\over 2t}\Bigr)\biggr)<+\infty.
\end{equation}

Assume that the curvature tensor satisfies the following growth condition
\begin{equation}\label{growth}
|||\Omega_u|||+\sum_{i=1}^d |||(L_{H_i}\Omega)_u|||\leq C\, \bigl(1+ \rho(x, \pi(u))^2\bigr)
\end{equation}
where $L_{H_i}$ denotes the Lie derivative with respect to $H_i$.
\vskip 2mm

Let $v\in\H$, consider the functional $F_T: W_x^T(M)\ra \R$ defined by

\begin{equation*}\label{eq4.1}
F_T(\gamma(w))=\int_0^T\langle \dot v(t), dw_t\rangle.
\end{equation*}

Let $h\in\H$; then by \eqref{A3}, we have (see also \cite{Fa2})
\begin{equation}\label{eq4.2}
(D_hF_T)(\gamma(w))=\int_0^T\langle \dot v(t), q(t, h)dw_t\rangle + \int_0^T \langle \dot v(t), \dot{\hat h}_t(w)\rangle\, dt.
\end{equation}

Let $a\in \R^d$ and consider $v(t)=t a$ with $|a|=1$ in \eqref{eq4.2}, we have
\begin{equation}\label{eq4.3}
(D_hF_T)(\gamma(w))=-\int_0^T\langle q(t, h)a, dw_t\rangle + \int_0^T \langle a, \dot{\hat h}_t(w)\rangle\, dt.
\end{equation}

Let $\{e_1, \cdots, e_d\}$ be an orthonormal basis of $\R^d$; define
$$\dis C_i(w,t,\tau)=-\int_\tau^t \Omega_{u_s(w)}\bigl(e_i, \circ dw(s)\bigr)\, {\bf 1}_{(\tau<t)}. $$

Then by Fubini theorem, the term $q(t,h)$ has the expression
\begin{equation*}
q(t,h)=-\sum_{i=1}^d \int_0^T \dot h^i(\tau) C_i(w,t,\tau)\, d\tau.
\end{equation*}

According to \eqref{eq4.3}, the gradient $D_\tau F_T$ has the following expression:

\begin{equation}\label{eq4.4}
(D_\tau F_T)(\gamma(w))=\sum_{i=1}^d \Bigl(\int_\tau^T \langle C_i(w,s,\tau)a, dw_s\rangle\Bigr) e_i + a + {1\over 2}\int_\tau^T \ric_Z(u_s)\,a\,ds.
\end{equation}

We have
\begin{equation}\label{eq4.5}
\hbox{\rm Var}(F_T)=\mathbb{E}(F_T^2)-\mathbb{E}(F_T)^2=|a|^2T=T.
\end{equation}

\begin{proposition}\label{prop4.1} Assume \eqref{growth}. Let
\begin{equation*}
\chi_T={\mathbb{E}\Bigl(\int_0^T |D_\tau F|^2\, d\tau\Bigr)\over \hbox{\rm Var}(F_T)}.
\end{equation*}
Then
\begin{equation}\label{eq4.6}
\chi_T= 1+{T\over 2}\langle \ric_Z(u_0)a,a\rangle+o(T)\quad\hbox{as } T\ra 0
\end{equation}
where $u_0$ is the initial point of \eqref{SDE2}.
\end{proposition}

\vskip 2mm
{\bf Proof.} We have, using \eqref{eq4.4},
\begin{equation*}
\begin{split}
|D_\tau F_T|^2&=\sum_{i=1}^d \Bigl(\int_\tau^T \langle C_i(w,s,\tau)a, dw_s\rangle\Bigr)^2 + |a|^2 +{1\over 4}\Bigl|\int_\tau^T \ric(u_s)a\, ds\Bigr|^2\\
&+\Big\langle a, \int_\tau^T \ric(u_s)a\,ds\Big\rangle + 2\sum_{i=1}^d \int_\tau^T \langle C_i(w,s,\tau)a, dw_s\rangle\, a^i\\
&+2\int_0^d \int_\tau^T \langle C_i(w,s,\tau)a, dw_s\rangle\cdot \int_\tau^T \langle \ric(u_s)a, e_i\rangle\, ds.
\end{split}
\end{equation*}

Put respectively
\begin{equation*}
{\mathbb E}\Bigl(\int_0^T |D_\tau F_T|^2\, d\tau\Bigr)=I_1(T)+I_2(T)+I_3(T)+I_4(T)+I_5(T)+I_6(T).
\end{equation*}

It is obvious that $I_2(T)=|a|^2 T=T$ and $I_5(T)=0$. We have
\begin{equation*}
I_1(T)=\sum_{i=1}^d \int_0^T\Bigl(\int_\tau^T {\mathbb E}(|C_i(w,s,\tau)a|^2)\, ds\Bigr)\, d\tau.
\end{equation*}

Now by growth condition \eqref{growth} and \eqref{eq4.0}, there is a constant $\delta>0$ such that
\begin{equation}\label{eq4.7}
 {\mathbb E}(|C_i(w,s,\tau)a|^2)\leq \delta\, (s-\tau).
 \end{equation}
So that $I_1(T)\leq \delta T^3/6$. By condition \eqref{eq2.02}, it is easy to see that $I_3(T)\leq {K_1^2T^3\over 12}$. It follows that
$I_6(T)\leq {\sqrt{\delta} K_1\over 6}T^3$. Now for $I_4(T)$, we have

\begin{equation*}
\lim_{T\ra 0} {I_4(T)\over T^2}= {1\over 2}\langle \ric (u_0)a,a\rangle.
\end{equation*}

Combining these estimates together with \eqref{eq4.5}, we get \eqref{eq4.6}. $\square$

\vskip 2mm
\begin{thm}\label{th4.1} Assume \eqref{eq2.02} and \eqref{growth}. Let $K_2(x)$ be the lower bound of $\Ric_x$. Then as $T\ra 0$,
\begin{equation}\label{eq4.8}
1-{K_1T\over 2} + o(T)\leq Spect(\L)\leq 1+{K_2(x)T\over 2} + o(T).
\end{equation}
\end{thm}

\vskip 2mm
{\bf Proof.}
For $K_2>0$, set $\beta={K_1\over K_2}$. As $T\ra 0$, we have

$$\aligned&\sqrt{(2+\beta)\Big(2+2\beta-\beta\e^{-\frac{K_2T}{2}}\Big)}=\sqrt{(2+\beta)^2\Big(1+\frac{\beta}{2+\beta}\Big(1-\e^{-\frac{K_2T}{2}}\Big)\Big)}\\
&=(2+\beta)\sqrt{1+\frac{\beta}{2+\beta}\frac{K_2T}{2}+o(T)}\\
&=(2+\beta)\Big(1+\frac{\beta}{2+\beta}\frac{K_2T}{4}+o(T)\Big).
\endaligned$$
So, for $K_2>0$, as $T\ra 0$,
$$\aligned&\beta\sqrt{(2+\beta)\Big(2+2\beta-\beta\e^{-\frac{K_2T}{2}}\Big)}\ \e^{-\frac{K_2T}{4}}\\
&=\beta(2+\beta)\Big(1+\frac{\beta}{2+\beta}\frac{K_2T}{4}+o(T)\Big)\Big(1-\frac{K_2}{4}T+o(T)\Big)\\
&=\beta(2+\beta)\Big[1+\frac{T}{4}\Big(\frac{K_1}{2+\beta}-K_2\Big)+o(T)\Big]
\\&=\beta(2+\beta)\Big[1-\frac{K_2 T}{2(2+\beta)}+o(T)\Big].
\endaligned$$

By (\ref{c3.2}), we get
$$Spect(\L)^{-1}\leq(1+\beta)^2-\beta(2+\beta)\Big[1-\frac{K_2 T}{2(2+\beta)}+o(T)\Big]=1+\frac{K_1 T}{2}+o(T),$$
which implies that
$$Spect(\L)\geq 1-\frac{K_1 T}{2}+o(T).$$

For $K_2<0$, by (\ref{c3.22}),
$$\aligned Spect(\L)^{-1}&\leq\frac{1}{2}+\frac{1}{2}\bigg(1+K_1\frac{1-\e^{-\frac{K_2T}{2}}}{K_2}\bigg)^2
=\frac{1}{2}+\frac{1}{2}\Big(1+\frac{K_1}{K_2}\Big(\frac{K_2 T}{2}+o(T)\Big)\Big)^2\\
&=1+\frac{K_1 T}{2}+o(T),\endaligned$$
which implies again
$$Spect(\L)\geq 1-\frac{K_1 T}{2}+o(T).$$

Now in \eqref{eq4.6}, taking the vector $a$ such that $\ric(u_0)a=K_2(x) a$ yields
$$\dis Spect(\L)\leq 1 + \frac{K_2(x) T}{2} + o(T).$$

The proof of \eqref{eq4.8} is completed. $\square$

\begin{cor} Assume \eqref{growth}. In the case where $\Ric =-K_1\Id$ with $K_1\geq 0$, we have
\begin{equation*}
\Bigl|Spect(\L)-1+{K_1 T\over 2}\Bigr|=o(T) \quad \hbox{\rm as }T\ra 0.
\end{equation*}
\end{cor}

\beg{thebibliography}{99}

\leftskip=-2mm
\parskip=-1mm

\bibitem{Aida} S. Aida, Gradient estimates of harmonic functions and the asymptotics of spectral gaps on path spaces, {\it Interdisplinary Information Sciences}, 2 (1996), 75-84.

\bibitem{AidaElworthy} S. Aida and D. Elworthy, Differential calculus on path and loop spaces I. logarithmic Sobolev inequalities on path spaces,
{\it C. R. Acad. Sci. Paris}, 321 (1995), 97-102.

\bibitem{BakryEmery} D. Bakry and M. Emery, Diffusion hypercontractivities, {\it S\'em. de Probab.}, XIX, Lect. Notes in Math., 1123 (1985), 177-206,
Springer.

\bibitem{BakryLedoux} D. Bakry and M. Ledoux, A logarithmic Sobolev form of the Li-Yau parabolic inequality,
{\it Rev. Mat. Iberoamericana}, 22 (2006), 683-702.

\bibitem{Bismut} J. M. Bismut,  \emph{Large deviation and Malliavin Calculus},
Birkh\"auser, Boston/Basel, 1984.


\bibitem{CHL} B. Capitaine,  E. P. Hsu and M. Ledoux, \emph{Martingale representation and a simple proof of logarithmic
Sobolev inequalities on path spaces,} Elect. Comm. Probab. 2(1997),
71--81.

\bibitem{CM} A. B. Cruzeiro and P. Malliavin, \emph{ Renormalized Differential
Geometry on path space: Structural equation, Curvature}.  J.
Funct. Anal. 139 (1996), p.119-181.

\bibitem{CF} A. B. Cruzeiro and S. Fang, \emph{Weak Levi-Civita connection for the
damped metric on the Riemannian path space and Vanishing of Ricci tensor in adapted differential
geometry}, J. Funct. Anal. 185 (2001), 681-698.

\bibitem{Dr} B. Driver,  \emph{A Cameron-Martin type quasi-invariant theorem
for Brownian motion on a compact Riemannian manifold},  J. Funct.
Anal. 110 (1992), 272--376.


\bibitem{ElworthyLi} D. Elworthy and X.M. Li, Formulae for the derivatives of heat
semi-group. {\it J. Funct. Anal.} 125 (1994), 252--287.

\bibitem{ELJL} K.D. Elworthy, Y. Le Jan and X.M. Li, \emph{on the geometry
of diffusion operators and stochastic flow}, Lect. notes in Math.,
1720 (1999), Springer.

\bibitem{ElworthyYan} K.D. Elworthy and Y. Yan, The Vanishing of harmonic one-forms on base path spaces,
{\it J. Funct. Analysis}, 264 (2013, 1168-1196.

\bibitem{ES} O. Enchev and D. Stroock, \emph{Towards a Riemannian geometry on
the path space over a Riemannian manifold}, J. Funct. Anal. 134
(1995), p. 392-416.


\bibitem{Fang1} S. Fang, \emph{In\'egalit\'e du type de Poincar\'e sur l'espace
des chemins riemanniens}, C.R. Acad. Sci. Paris, 318 (1994),
257-260.

\bibitem{Fa2} S. Fang, \emph{Stochastic anticipative integrals on a
Riemannian manifold}, J. Funct. Anal. 131 (1995), 228-253.

\bibitem{FangMalliavin}  S. Fang and P. Malliavin, Stochastic analysis on the path space of a Riemannian manifold, {J. Funct. Analysis}, 131 (1993), 249-274.

129(2005), 339--355.

\bibitem{FWW} S. Fang, F.Y. Wang and B. Wu, Transportation-cost inequality on path spaces with uniform distance, {\it Stochastic Process. Appl.} 118 (2008), no. 12, 2181¨C2197.


\bibitem{Hs1} E. P. Hsu, \emph{Quasi-invariance of the Wiener measure on the path
space over a compact Riemannian manifold}, J. Funct. Anal. 134
(1995), p. 417-450.

\bibitem{Hs2} E. P. Hsu, \emph{Quasi-invariance of the Wiener measure on path
spaces: Noncompact case}, J. Funct. Anal. 193 (2002), p. 278-290.

\bibitem{Hs3} E. P. Hsu, \emph{Logarithmic Sobolev inequalities on
path spaces over Riemannian manifolds, } Comm. Math. Phy. 189
(1997), 9-16.



\bibitem{Naber} A. Naber, Characterizations of bounded Ricci curvature on smooth and nonsmooth spaces, {\it arXive: 1306.651}.

\bibitem{Stroock} D. Stroock, {\it An introduction to the analysis of paths on a
Riemannian manifold}, Math. Survey and Monographs, vol. 74, AMS,
2000.

\end{thebibliography}

\end{document}